\documentclass[SEDP]{cedram-sem}

\usepackage{amsfonts}
\usepackage{amsthm}
\theoremstyle{plain}

\newtheorem{thm}{Theorem} 

\theoremstyle{definition}

\newtheorem{df}{Definition} 

\usepackage{cite}
\usepackage{comment}
\usepackage{amsmath,bm}
\usepackage{amscd}
\usepackage{amssymb}
\usepackage[pdftex]{color,graphicx,hyperref}
\usepackage{latexsym}
\usepackage{color}
\usepackage{cases}
\usepackage{xfrac}
\usepackage[normalem]{ulem}
\usepackage{environ}
\usepackage{xcolor}
\usepackage{mathtools}
\usepackage{stackengine}
\usepackage[framemethod=tikz]{mdframed}
\mdfdefinestyle{MyFrame}{%
	linecolor=black,
	outerlinewidth=1pt,
	roundcorner=5pt,
	innertopmargin=\baselineskip,
	innerbottommargin=\baselineskip,
	innerrightmargin=10pt,
	innerleftmargin=10pt,
	backgroundcolor=white!20!white}

\usepackage{tikz}
\usetikzlibrary{positioning}

\newcommand{\CX}{{\mathcal X}}
\newcommand{\CY}{{\mathcal Y}}
\newcommand{\CV}{{\mathcal V}}
\newcommand{\intO}[1]{\int_{\mathbb{T}^d} #1 \ \dx}

\newcommand{\dx}{{\rm d} {x}}
\newcommand{\dy}{{\rm d} {y}}
\newcommand{\dt}{{\rm d} t }

\newcommand{\vr}{\varrho}
\newcommand{\vrn}{\varrho_n}
\newcommand{\vw}{\vc{w}}
\newcommand{\vu}{\vc{u}}
\newcommand{\vv}{\vc{v}}
\newcommand{\pt}{\partial_{t}}
\newcommand{\px}{\partial_x}
\newcommand{\Div}{\operatorname{div}}
\newcommand{\Grad}{\nabla}
\newcommand{\eq}[1]{\begin{equation}
\begin{split}
#1
\end{split}
\end{equation}}
\newcommand{\eqh}[1]{\begin{equation*}
\begin{split}
#1
\end{split}
\end{equation*}}
\newcommand{\ep}{\varepsilon}
\newcommand{\R}{\mathbb{R}}
\newcommand{\T}{\mathbb{T}}
\newcommand{\Td}{\mathbb{T}^d}
\newcommand{\lr}[1]{\left( #1 \right)}
\newcommand{\vc}[1]{{\bf{ #1}}}
\newcommand{\Nu}{\mathcal{V}}

\newcommand{\rtrw}{\widetilde{\sqrt{\varrho} \vc{w}}}
\newcommand{\lel}{\left\langle}
\newcommand{\ril}{\right \rangle}

\newcommand{\DQ}{\widetilde{D_Q} }
\newcommand{\levert}{\left\vert}
\newcommand{\rivert}{\right \vert}
\newcommand{\bvr}{\bar{\varrho}}
\newcommand{\bvw}{\bar{\vc{w}}}
\newcommand{\Dt}{\frac{ d}{dt}}

\title
[On the  dissipative Aw-Rascle system]
{Recent advances in the analysis of the  dissipative Aw-Rascle system}

\alttitle{Avancées récentes dans l'analyse du système dissipatif Aw-Rascle}

\author
{\firstname{Ewelina}  \lastname{Zatorska}}

\address{Mathematics Institute, \\
University of Warwick\\
Zeeman Building, \\
Coventry CV4 7AL,\\ 
United Kingdom}
\thanks{The work of the author was  supported by the EPSRC Early Career Fellowship no. EP/V000586/1.}

\email{ewelina.zatorska@warwick.ac.uk}
\keywords{Dissipative Aw-Rascle system, weak-solutions, Young measures, convex integration, duality solutions, hard congestion limit}
  
\altkeywords{Système dissipatif d'Aw-Rascle, solutions faibles, mesures de Young, intégration convexe, solutions de dualité, limite de congestion dure}

\subjclass{35L65, 35Q31, 76N10}

\begin{document}


\begin{abstract}
 The one-dimensional Aw-Rascle (AR) system has become a cornerstone of macroscopic models for single-lane vehicular traffic. A possible generalization of this model to a multi-dimensional setting is the so-called dissipative AR model, which is more suited to capturing crowd dynamics. This review summarizes recent studies that analyze the dissipative AR model, its hard congestion limit, the non-uniqueness of weak solutions, the existence and asymptotics of solutions within the duality framework, non-local interactions, and the existence of regular solutions.
\end{abstract}

\begin{altabstract}
  Le syst\`{e}me Aw-Rascle (AR) unidimensionnel est devenu un pilier des mod\`{e}les macroscopiques pour le trafic v{\'e}hiculaire \`{a} une seule voie. L'une des g{\'e}n{\'e}ralisations possibles de ce mod\`{e}le dans un cadre multidimensionnel est le mod\`{e}le AR dissipatif, mieux adapt{\'e} pour d{\'e}crire la dynamique des foules. Cette revue r{\'e}sume les travaux r{\'e}cents qui analysent le mod\`{e}le AR dissipatif, sa limite de congestion forte, la non-unicit{\'e} des solutions faibles, l'existence et l'asymptotique des solutions dans le cadre de la dualit{\'e}, les interactions non-locales et l'existence de solutions r{\'e}guli\`{e}res.
\end{altabstract}

\maketitle


\section{Introduction}

The development of macroscopic models of vehicular traffic originates in the works of Lighthill and  Whitham \cite{LW} and Richards \cite{R} from the late 50s. These first-order models, taught nowadays in the undergraduate mathematics courses on transport equation, were basically equations for conservation of mass of the vehicles with the velocity of motion given explicitly as a function of the density. Their second-order correction, with additional equation for velocity, appeared more than decade later in the work of Payne \cite{Payne} and Whitham \cite{Whitham}. 
The subsequent history of these fluid-like models of traffic experienced a few plot twists. In 1995  their ``requiem''  was announced by Daganzo \cite{Daganzo}, who based his criticism on the three differences between the traffic and the fluid particles: the anisotropy of interactions, the finite extent of the shock waves, and personal preferences of the drivers. 
These shortcomings  were subsequently addressed by Aw and Rascle   \cite{AR} and independently by Zhang \cite{Zhang}. They ``resurrected'' the fluid-like traffic models by including in the system certain function describing the anticipation of road conditions in front of the drivers $P$. In short, the introduced model, called the Aw-Rascle (AR) model, resembled the compressible pressureless Euler equations, except that  the conservation of momentum involved not the actual velocity of motion $u$, but the preferred velocity $w$
\begin{subnumcases}{\label{AR}}
\pt \vr+\px(\vr u)=0,\label{AR1}\\
\pt(\vr w)+\px(\vr uw)=0,\label{AR2}\\
w=u+P(\vr)
\label{AR3}
\end{subnumcases} 
where $\vr(t,x)$ stands for the density, i.e.  a number of vehicles per unit length of road.  The anticipation function $P=P(\vr)>0$ is the offset function between the velocities, and it expresses the observation  that the actual velocity is always smaller than the  desired one due to congestion of the cars ahead. 
The one-dimensional AR system has been derived in \cite{AwKlar} from the  Follow-the-Leader particle model with a specific form of the offset function $P(\vr)=\vr^\gamma$.
Later on,  in \cite{BDDR}, a singular offset function was considered 
\eq{\label{sing_cost}
P(\vr)=\ep\lr{\frac{1}{\vr}-\frac{1}{\bar\vr}}^{-\gamma}}
 with the maximal density constraint $\bar\vr>0$.
This form of the offset function causes the density $\vr$ to stay always below its critical value $\bar\vr$, provided that it was so initially. As demonstrated in \cite{BDDR} the limit $\ep\to 0$ in the system \eqref{AR} leads to a constrained pressureles gas dynamics system. Similarly to the pressure in the fluid mechanic equations, $P(\vr)$ propagates the flow perturbations. Compressible Euler and Navier-Stokes equations with singular pressure functions \eqref{sing_cost}, and the corresponding hard-congestion limits were considered, for example, in \cite{PZ, DMZ, Perrin}.

In several space dimensions, the offset function is no longer a scalar. In the two-dimensional model of traffic with driving direction analyzed in \cite{HMV}, the a vector of offset functions in each of directions was considered $\vc{P}(\vr)=[P_1(\vr),P_2(\vr)]$. Another multidimensional generalization of the AR model was proposed in \cite{ABDM} to describe the motion of pedestrians. In this model the scalar offset function $P(\vr)$ was replaced by a gradient of some cost function $\Grad p(\vr)$,
so that the velocity offset is also a vector, i.e. $\vw=\vu+\Grad p(\vr)$. We then arrive at the following multi-dimensional version of \eqref{AR}:
\begin{mdframed}[style=MyFrame]
\begin{subnumcases}{\label{AR_multi}}
\pt \vr+\Div (\vr \vu)=0,\label{AR1_multi}\\
\pt(\vr \vw)+\Div(\vr \vu\otimes\vw)=0,
\label{AR2_multi}\\
\vw=\vu+\Grad p(\vr).\label{AR3_multi}
\end{subnumcases} 
	\end{mdframed}
The research content of \cite{CGZ,CMPZ, CNPZ, CPSZ, CPZ, CFZ, M}, that we summarize in the current overview, concerns analysis of this system, called the {\emph{dissipative AR system}}. To see where the dissipation appears, let us rewrite the momentum equation in terms of the unknowns $\vr$ and $\vu$, we then (formally) obtain
\begin{align}\label{NS-multi}
		\partial_t (\vr \vc{u}) + \Div (\vr \vc{u} \otimes \vu ) = \nabla_x (\vr Q^\prime (\vr)\Div \vu) + \mathcal{L}[\nabla_x Q(\vr), \nabla_x \vu] ,
\end{align}
where $ Q^\prime(\vr) = \vr p^\prime (\vr)$ and $$\mathcal{L}[\nabla_x Q(\vr), \nabla_x \vu] = \nabla_x(\nabla_x Q(\vr) \cdot \vu )- \Div(\nabla_x Q(\vr) \otimes \vu ),$$ 
which is a lower order term. Indeed,  a simple calculation yields
\begin{align*}
     \left(\mathcal{L}[\nabla_x Q(\vr), \nabla_x \vu]\right)_j= \sum_{i=1}^{3} \left( \partial_{x_i}Q(\vr) \partial_{x_j}u_{i} - \partial_{x_j} Q(\vr) \partial_{x_i} u_i   \right)  \text{ for } j=1,2,3.
\end{align*}
In particular, in one-dimensional setting, equation \eqref{NS-multi} reduces to
\eq{\label{NS_1D}
\partial_t(\vr u)+\partial_x(\vr u^2)-\partial_x\lr{\vr^2 p'(\vr)\partial_x u}=0,}
which is just the momentum equation of the compressible, pressureless Navier-Stokes equations with density-dependent viscosity coefficient $\mu(\vr)=\vr^2p'(\vr)$. Similar equation has been extensively studied in the one and multi-dimensional setting, see for example \cite{Haspot, HZ} and the references therein.

The paper is organised as follows. In the  first part, we discuss the results on  the multi-dimensional dissipative AR system \eqref{AR_multi} with the offset function equal to $\Grad \vr^\gamma$ with fixed  $\gamma\geq 1$. Here we report on four main results: the short-time existence of regular solutions, the global-in-time existence of measure-valued weak solutions, the weak(measure-valued)-strong uniqueness of solutions, and the ill-posedness of the dissipative AR system in the class of weak solutions, form the papers \cite{CGZ, CPZ, CFZ}. 

Second part of this review is dedicated to discussion of connections between the dissipative AR system and other systems with two-velocity structure. For example, we discuss the version of the dissipative AR model, in which the offset function is non-local and depends on the macroscopic density through convolution with the offset kernel $K(x)$, i.e. $\vu=\vw-\Grad K \ast  \vr$, analysed in \cite{CPSZ}. Another example is the dissipative AR system is written in one-dimensional form. For a suitable singular cost function $p(\vr)$, the system coincides with a moedel of lubrication derived in \cite{LM}
by Lefebvre-Lepot and B. Maury. We discuss the existence of regular solutions for $\ep$ fixed, and the hard congestion limit passage within the regime of weak and duality solutions, following the results of \cite{CNPZ,CMPZ,M}.

\section{Dissipative AR system in several space-dimensions}
\subsection{Existence of regular solutions}
The first question that we address in this note is the existence of sufficiently regular solution that satisfy the system \eqref{AR_multi} in the strong sense.
For simplicity, we consider system \eqref{AR_multi} on the three-dimensional torus $\T^3$,
with the initial data
\begin{equation*}
    \vr(0,x)=\vr_0(x),\qquad \vu(0,x)=\vu_0(x),
\end{equation*}
and with the cost function of the form
$p(\vr)=\vr^\gamma$ with $\gamma\geq 1$.

For $T>0$ and $k \in \mathbb{N}$ we introduce the notation: 
\begin{equation*} \label{def:CX}
\begin{aligned}
&\CX_k(T):=L^2(0,T;H^k(\T^3)) \cap L_\infty(0,T;H^{k-1}(\T^3)),\\
& \CY_k(T):=\{ f\in L^\infty(0,T;H^k(\T^3)): \; \pt f \in L^\infty(0,T;H^{k-1}(\T^3))\\
& \CV_k(T):=\{ g \in L^2(0,T;H^{k+1}(\T^3)) \cap L^\infty(0,T;H^k(\T^3)): \;\pt g \in L^2(0,T;H^{k-1}(\T^3)) \}
\end{aligned}
\end{equation*}
with norms defined in a natural way as appropriate sums of norms.

The following existence result is the main result of the joint work with Chaudhuri and Piasecki,  \cite{CPZ}.
\begin{thm} \label{thm:main}
Assume the initial data satisfies $\vr_0>0$, and either
\begin{equation} \label{init1}
(\vr_0,\vu_0) \in H^4(\T^3) \times H^3(\T^3)
\end{equation}
or
\begin{equation} \label{init2}
\vr_0 \in H^3(\T^3), \quad \vu_0+\nabla p(\vr_0) \in H^3(\T^3).    
\end{equation}

Then there exists $T>0$ such that system \eqref{AR_multi} admits a unique solution $(\vr,\vw) \in \CV_3(T)\times \CY_3(T)$ with the estimate
$$
\|\vr\|_{\CV_3(T)}+\|\vw\|_{\CY_3(T)} \leq C(\|\vr_0\|_{H^4(\T^3)},\|\vu_0\|_{H^3(\T^3)})
$$
in case of \eqref{init1} or 
$$
\|\vr\|_{\CV_3(T)}+\|\vw\|_{\CY_3(T)} \leq C(\|\vr_0\|_{H^3(\T^3)},\|\vu_0+\nabla p(\vr_0)\|_{H^3(\T^3)})
$$
in case of \eqref{init2}.
\end{thm}

The  solutions to  \eqref{AR_multi} are first approximated by solutions to the following iterative scheme 
\begin{equation*} \label{iter}
\left\{ \begin{array}{lr}
\vr^{n+1}_t + \Div(\vr^{n+1}\vw^{n+1})=\Div(\vr^n p'( \vr^n)\nabla\vr^{n+1}),\\
\vw^{n+1}_t + \vu^n\cdot \nabla \vw^{n+1} = 0,\\
(\vr^{n+1},\vw^{n+1})|_{t=0}= (\vr_0,\vu_0+\nabla p(\vr_0)).
\end{array}\right.
\end{equation*}
Given $(\vr^n,\vw^n)$, one defines $\vu^n=\vw^n-\nabla p(\vr^n)$ and solves the second equation of
\eqref{iter} for $\vw^{n+1}$. Next, one may use $\vw^{n+1}$ to determine the solution $\vr^{n+1}$ to the first equation. The problem is therefore reduced to solving two separate liner equations  at each iteration. Convergence of this iterative scheme is then proved using the Banach fixed point theorem.

The key component of the proof are  the $L^p$ estimates for the transport equations involved. The departure point to obtain them  is the explicit solution formula following from the method of characteristics. Partial results of this type have been used in the theory of compressible Navier-Stokes equations (see \cite{KNP} and references therein). We extend these results to systems of mixed hyperbolic-parabolic type, which exhibit dissipation in the continuity equation but lack it in the momentum equation. Indeed, dissipativity in the first equation of system \eqref{iter} gives parabolic estimates, but a delicate part is to ensure positivity of the solution at each step of the iteration. 

\subsection{Global-in-time solutions}
The previous result was restricted to the regular solutions on possibly small time interval. In order to extend these solutions to large times, one would need to make a restriction on the size of initial data.  For long-time solutions with arbitrary large initial data, one usually needs to look for solution in a bigger class of functions risking loss of regularity or uniqueness of solution. In the case of dissipative AR system it turns out that the usual class of distributional weak solutions is still to narrow to show compactness of approximating sequences satisfying the energy-type a-priori estimates.
To motivate the definition of wider class of solutions, the measure-valued solution, let us first observe that at the level of sufficiently regular solutions, and under constraint $\vw=\vu+\Grad p(\vr)$ the following reformulation
\begin{subnumcases}{\label{sym}}
		\pt \vr+\Div(\vr \vw)-\Div(\sqrt{\vr}\Grad Q)=0, \\
		\pt(\vr \vw)+\Div(\vr \vw\otimes\vw)=\Div(\Grad Q\otimes\sqrt\vr\vw),
\end{subnumcases}
{where} $Q'(\vr)=\sqrt{\vr}p'(\vr)$ are equivalent to system \eqref{AR_multi}. In particular, the first equation \eqref{sym} resembles the porous medium equation. It is therefore not difficult to deduce the following global-in-time a-priori estimates for a sequence of regular solutions to either of the systems \eqref{AR_multi} or \eqref{sym}
\begin{align}
&\|\vr_n\|_{L^\infty(0,T; L^1(\T^d))}\leq C,\label{est_mas}\\
&\|\vr_n|\vw_n|^2\|_{L^\infty(0,T; L^1(\T^d))}\leq C,\label{est_mom}\\
   & \|E(\vr_n)\|_{L^\infty(0,T; L^1(\T^d))}\leq C, \label{eng:1}\\
   & \|\Grad Q(\vr_n)\|_{L^2((0,T)\times\T^d)}\leq C, \label{eng:2}
\end{align}
where  $d$ stands for a spacial dimension and it is equal 2 or 3.
In the estimates above, $E$ denotes the energy associated with the porous medium structure, i.e.
\eqh{\label{def:E}
E(\vr)=\int_0^\vr p(s)\ {\rm d}s,}
and $C$ is a constant  that  only on the initial data. More precisely, we require the initial data to be energy-bounded
\eqh{\label{MV:IE:2}
        &0\leq  \varrho_0\quad \text{ in } \T^d\quad \text{ and }\quad \intO{\left( \frac{1}{2} \vr_0 |\vw_0|^2+ E(\vr_0)\right)} < \infty .}
 The proof of the existence of solutions to the system \eqref{AR_multi}  satisfying these a-priori estimates consists of a multi-layer construction procedure, involving, and the basic level, a version of Brenner's model introduced and studied in  \cite{B1,B2,FV}. 
 The approximation procedure eventually gives rise to a sequence of solutions $\vr_n,\vu_n,\vw_n$ that generate a generalized measure-valued solution. Indeed, the a-priori estimates (\ref{est_mas}-\ref{eng:2}) lead to the following estimates of consecutive terms in the distributional formulations of the equations:
 \begin{itemize}
\item {in the continuity equation:}
\eqh{
\pt \vr_n+\Div(\underbrace{\sqrt{\vr_n}\sqrt{\vr_n} \vw_n}_{L^\infty(L^p)})-\Div(\underbrace{\sqrt{\vr_n}\Grad Q(\vrn)}_{L^2(L^p)})&=0,
}
\item {in the momentum equation:}
\eqh{\pt(\sqrt{\vr_n}\sqrt{\vr_n} \vw_n)+\Div(\underbrace{\sqrt{\vr_n}\vw_n\otimes\sqrt{\vr_n}\vw_n}_{L^\infty(L^1)})&=\Div( \underbrace{\Grad Q(\vr_n)\otimes\sqrt{\vr_n}\vw_n}_{L^2(L^1)}).}
\end{itemize}

In particular, the nonlinear terms in the momentum equation are bounded only in $L^1$, and so, their limits may end up to be measures only. This problem motivates the use of Young measures which are tailored to identify the limits in the nonlinear terms. The situation is quite analogous  to the compressible Euler system considered in \cite{Basaric}, where the  measure-valued solution was generated by the sequence $ \{\vr_n, \vr_n \vu_n\}_{n\in \mathbb{N}}$. In our case, however,  there is an additional nonlinear term $ \Grad Q(\vrn)\otimes\sqrt{\vr_n}\vw_n $ which is only a-priori bounded in $L^1$ and its limit cannot be identified in this space.  The solution to \eqref{AR_multi} is therefore a Young measure $ \Nu $, generated by $ \left\{(\vr_n, \sqrt{\vr_n} \vw_n, \Grad Q(\vr_n))\right\}_{n\in \mathbb{N}} $, where $(\vr_n, \vw_n)$ is the solution of a suitable approximation of \eqref{sym}. More precisely, we have the following definition.
\begin{df}\label{MVdef}
		We say that a parametrized measure $\{ \Nu_{t,x} \}_{(t,x)\in (0,T)\times \T^d}$,
	\begin{align*}
		\Nu \in L^{\infty}_{\text{weak-(*)}} \big( (0,T)\times \T^d ;\mathcal{P}(\mathcal{F} )\big),
	\end{align*}
	on the phase space
	\begin{align*}
		\mathcal{F}= \left\{\left(s, \rtrw, \mathbf{F}\right) \mid s\in [0,\infty),\; \rtrw \in \R^d,\; \vc{F} \in \R^d \right\}
	\end{align*}
	is a measure--valued solution of the dissipative  AR system \eqref{sym} in $(0,T)\times \T^d$, with the bounded energy initial data \eqref{MV:IE:2} and dissipation defect $\mathcal{D}$,
	\begin{align*}
		\mathcal{D}\in L^{\infty}(0,T),\; \mathcal{D}\geq 0,
	\end{align*}
	if the following conditions hold:
	\begin{itemize}
		\item The unknowns enjoy the following regularity: 
		\begin{align*}
			&  \vr=\lel \Nu_{t,x}; s  \ril  \in C_{\text{weak}}(0,T;L^{\gamma+1}(\T^d)),\\
			& 	{\vr}^\alpha=\lel \Nu_{t,x}; {s}^\alpha  \ril  \in C_{\text{weak}}(0,T;L^{\frac{1}{\alpha}(\gamma+1)}(\T^d))  \text{ with } 0<\alpha< \gamma+1 ,\\
			&\lel \Nu_{t,x};  \rtrw \ril \in L^\infty(0,T;L^2(\T^d)) \text{ with } \sqrt{\vr} \lel \Nu_{t,x};  \rtrw \ril =  \lel \Nu_{t,x}; \sqrt{s} \rtrw \ril,\\
			&Q(\vr)=\lel \Nu_{t,x}; Q(s)  \ril \in L^2(0,T; W^{1,2}(\T^d)) \text{ with } \Grad Q(\vr)=\lel \Nu_{t,x}; \DQ  \ril.
		\end{align*}	
		\item For a.e. $\tau \in (0,T) $ and $\psi \in C^{1}([0,T]\times {\T^d})$ we have
		\begin{align*} 
			\begin{split}
				&\int_{\T^d}\vr(\tau, \cdot) \psi(\tau, \cdot) \dx - \int_{ \T^d} \vr_0 \psi(0, \cdot) \dx \\
				&\quad = \int_{0}^{\tau} \int_{ \T^d} \left[ \vr  \partial_{t}\psi + \sqrt{\vr}\lel \Nu_{t,x};  \rtrw \ril \cdot \Grad \psi - \sqrt{\vr} \Grad Q(\vr) \cdot \Grad \psi \right] \dx \dt.
			\end{split}
		\end{align*}
		\item There exists a defect measure $r^{M} \in L^\infty_{\text{weak-(*)}}(0,T;\mathcal{M}({\T^d};\R^{d\times d})) + \mathcal{M}([0,T]\times \T^d;\R^{d\times d})$ and $\xi\in L^{1}(0,T)$ such that for a.e. $\tau \in (0,T)$, for all $ \epsilon>0 $ and $\pmb{\phi}\in C^{1}([0,T]\times {\T^d};\R^d)$ we have	
		\begin{align*}\label{mv eqn 2}
			\begin{split}
				&\int_{ \T^d} \sqrt{\vr}(\tau,x) \lel\Nu_{\tau,x}; \rtrw  \ril \cdot \pmb{\phi}(\tau,\cdot) \dx - \int_{ \T^d} \lel \Nu_{0} ;  \sqrt{s} \rtrw  \ril \cdot \pmb{\phi}(0,\cdot) \dx\\
				&= \int_{0}^{\tau} \int_{ \T^d} \sqrt{\vr}  \left[ \lel \Nu_{t,x};  \rtrw   \ril \cdot \partial_{t} \pmb{\phi} +  \lel \Nu_{t,x}; \rtrw \otimes  \rtrw  \ril : \Grad \pmb{\phi}- \lel \Nu_{t,x}; \rtrw \otimes \DQ \ril : \Grad \pmb{\phi} \right] \dx \dt\\
				&  + \lel r^M;\Grad \pmb{\phi}\ril_{\mathcal{M}([0,\tau]\times \T^d), C([0,\tau]\times \T^d)} 
			\end{split}
		\end{align*}
	and
		\begin{equation*}\label{mom-def}
		\left\vert \lel r^M;\Grad \pmb{\phi}\ril_{\mathcal{M}([0,\tau]\times \T^d), C([0,\tau]\times \T^d)} \right\vert \leq  \left( 1+\frac{1}{4 \epsilon} \right) \int_0^\tau \xi(t) \mathcal{D}(t) \Vert \pmb{\phi} \Vert_{C^1({\T^d})} \dt + \epsilon \Vert \pmb{\phi} \Vert_{C^1([0,\tau]\times {\T^d})} \mathcal{D}(\tau).
	\end{equation*} 
		\item There exists a defect measure $ \mathcal{R} \in L^\infty_{\text{weak-(*)}}(0,T;\mathcal{M}(\T^d))+ \mathcal{M}([0,T]\times \T^d)$ such that for a.e. $ \tau \in (0,T) $ and $ \epsilon >0 $ we have
		\begin{align*}
			\begin{split}
				&\int_{ \T^d} \lel \Nu_{\tau ,x}; \frac{1}{2}  \levert \rtrw \rivert^2 + E(s) \ril \dx + \int_{0}^{\tau} 	\intO{ \lel \Nu_{t,x}; \vert \DQ \vert^2 \ril} \dt + \mathcal{D}(\tau) \\
				&\leq \int_{ \T^d} \lel \Nu_{0,x};  \frac{1}{2} \levert \rtrw \rivert^2 +E(s) \ril \dx + \int_{0}^{\tau} \intO{ \lel \Nu_{t,x}; \rtrw \cdot \DQ \ril }  \dt +\int_{(0,\tau)\times \T^d}{\text{d} \mathcal{R}}  ,
			\end{split}
		\end{align*}
		where
  \begin{align*}\label{En:def}
      \left\vert \int_{(0,\tau)\times \T^d}{\text{d} \mathcal{R}}   \right\vert \leq C \left(1+\frac{1}{4 \epsilon} \right) \int_{0}^\tau \mathcal{D}(t) \dt + \epsilon \mathcal{D}(\tau).
  \end{align*}
	Here we can identify $ \lel \Nu_{t,x}; E(s) \ril= E(\vr) $. 
	\end{itemize}
	\end{df}
	
As mentioned above, the existence of such solutions can be obtained by passing to the limit with subsequent parameters of approximation, and it has been described in \cite{CFZ}. here we recall only the following existence result
\begin{thm}\label{Theorem:1}
Let $p(\vr)=\vr^\gamma $, $ \gamma \geq1$, and let the initial data satisfy \eqref{MV:IE:2}. Then, for any $ T>0 $,  there exists a measure-valued solution to \eqref{sym} in the sense of 
 Definition \ref{MVdef}.
\end{thm}
Using dissipative effect on the density, we can in fact show that the projection of this Young measure onto the first component  is a Dirac delta. More precisely 
  \[ \Nu_{t,x}= \delta_{\{\vr(t,x)\}}\otimes Y_{t,x} \quad \text{for  a.a.\ } (t,x) \in (0,T) \times \T^d, \] 
    where $Y \in L^{\infty}_{\text{weak-(*)}} \big( (0,T)\times \T^d ;\mathcal{P}( \R^d \times \R^{d})\big)$. 
    We can also show that as long as a sufficiently regular solution exists, the measure-valued solution emanating from the same data is a Diract delta centered at the regular solution, i.e. the solutions coincide. This is the so-called weak-strong uniqueness result that can be formulated as follows
\begin{thm}\label{Theorem:2}
	Let  $ (\Nu, \mathcal{D}) $ be a measure-valued solution of \eqref{sym} in the sense of Definition \ref{MVdef}. Let $ (\bvr,\bvw) $ be a strong solution to the same system, s.t.
	\begin{equation*}
		\bvr \in C^1([0,T]; C^2(\T^d)),\; \bvw \in C^1([0,T]; C^2(\T^d; \R^d)) \quad
		\text{ with } \bvr >0,
		\end{equation*}
		 emanating from the initial data $ (\bvr_0,\bvw_0) \in (C^2(\T^d), C^2(\T^d;\R^d)) $ s.t. $ \bvr_0>0 $. If the initial data coincide, i.e.
\begin{equation*}
	\Nu_{0 ,x}= \delta_{\{ \bvr_0(x), \bvw_0(x)\}}, \text{ for a.e. }x\in \T^d
	\end{equation*} 
then $ \mathcal{D}=0 $, and 
\begin{align*}
		\Nu_{\tau ,x}= \delta_{\{ \bvr(\tau,x), \sqrt{\bvr}\bvw(\tau, x), \Grad Q(\bvr)(\tau,x) \}}, \text{ for a.e. }(\tau,x)\in(0,T)\times \T^d.
\end{align*}
\end{thm}
The proof of this theorem uses a generalization of the relative entropy method to measure the distance between the measure-valued solution defined above, and the regular solution from the previous section. For two smooth solutions $ (\bvr, \bar{\vw}) $ and $ (\vr,\vc{w}) $ this relative entropy functional has a form:
\eqh{\label{def:rel_ent}
\mathcal{E}(\varrho, \vc{w} \mid \bvr , \bar{\vw}) : = \intO{\left(\frac{1}{2} \varrho \vert \vw- \bar{\vw} \vert ^2 +  E(\varrho)-E(\bvr)-E^\prime (\bvr) (\varrho-\bvr) \right)},
}
but it needs to be modified if one of the solutions is of lower regularity. It is worth noticing that, unlike in the classical fluid systems, the first component measures the kinetic energy associated with the desired velocity of motion $\vw$, and not the actual one, while the potential part of the energy $E$ is not the  internal energy associated with the pressure, but with the porous medium diffusion in the continuity equation.

\subsection{Ill-posedness in the class of weak solutions}

The question of existence and uniqueness of weak solutions to system \eqref{AR_multi} was posed and answered in the joint work with Chaudhuri and Feireisl \cite{CFZ}. It was proven that any initial density--velocity data $(\vr_0, \vu_0) = (\vr(0, \cdot), \vu(0, \cdot))$ can connect to arbitrary terminal state $(\vr_T, \vu_T) = ((\vr(T, \cdot), \vu(T, \cdot))$ via a weak solution.
More specifically,  one requires regular  data compatible with mass and momentum conservation, i.e.
\begin{equation*} \label{cc1}
	\vr_0, \vr_T \in C^2(\T^d),\quad \inf_{\T^d} \vr_0 > 0,\quad \inf_{\T^d} \vr_T > 0, \quad \intO{ \vr_0 } = \intO{ \vr_T },
	\end{equation*}
together with 
\eqh{
	\vu_0, \vu_T \in C^3(\Td; \R^d),\quad	\intO{\vr_0 \vu_0 }= \intO{ \vr_T \vu_T } .
}
The drawback of this result is that the total energy of the system, defined here by
\[
E(\vr, \vu) = \frac{1}{2} \vr \left| \vu + \Grad p(\vr) \right|^2,
\] 
may experience jumps and actually increase in time. For sufficiently regular weak solutions this quantity should actually be conserved (because the linear momentum $\vr(\vu + \Grad p(\vr))$ is conserved), so one expects that the energy for the weak solutions would be at most dissipated. To obtain existence of arbitrary many weak solutions respecting the energy inequality, one needs to change the ansatz about the initial data. We have the following result.
\begin{thm}[{\bf Ill posedness in the class of admissible solutions}] \label{ccT2}
		Let $d=2,3$, and let  $p \in C^2((0,\infty))$.
		Let $\vr_0 \in C^2(\Td)$, $\inf_{\Td} \vr_0 > 0$ be given. 
		Then there exists an initial velocity $\vu_0 \in L^\infty(\Td; \R^d)$  such that
		the system \eqref{AR_multi}, endowed with the periodic boundary  conditions  admits infinitely many weak solutions in the class 
		\[
		\vr \in C^2([0,T] \times \Td), \vu \in L^\infty ((0,T) \times \Td; \R^d)
		\]
		satisfying
		\begin{equation*} \label{cc3b}
			\vr(0, \cdot) = \vr(T, \cdot) = \vr_0,\  
			(\vr \vu) (T, \cdot) = 0,
		\end{equation*}
		together with the energy inequality 		
\[		
	\Dt \intO{ E(\vr, \vu) } \leq 0,\ \intO{ E(\vr, \vu) (t, \cdot) } \leq 
\intO{ E(\vr_0, \vu_0 ) }.
\]	
	\end{thm}

	The proof of this theorem relies on application of the convex integration technique introduced by DeLellis and Sz\'ekelyhidi \cite{DeSz}, primarily to prove the existence of infinitely many {\emph{wild}} solutions to the incompressible Euler system. Subsequently, it was extended by Chiodaroli \cite{Ch} to the compressible Euler system, and more recently by  Buckmaster and Vicol to the incompressible Navier-Stokes equations \cite{BuVi}. It is not yet known if convex integration technique could be further extended to the compressible Navier-Stokes equations. Note, however, that weak inviscid limit of  compressible Navier--Stokes system with density-dependent viscosities has been recently used in \cite{CVY} to generate infinitely many global-in-time admissible weak solutions to the isentropic Euler system. Our application of the convex integration technique for dissipative system \eqref{AR_multi}, which in one-dimensional setting coincides with the compressible Navier-Stokes equations \eqref{NS_1D}, is thus an important step forward.

\section{Comparison with other two-velocity hydrodynamic models}
As mentioned in the previous section, the dissipative AR model might be viewed as an inviscid and pressureless version of the Brenner model, analysed by Feireisl and Vasseur in \cite{FV}.
The Brenner model postulates the existence of two velocities: $\vv_m$ -mass velocity and $\vv$-volume velocity. According to the author, in general $\vv_m\neq\vv$ and the inequality becomes the more significant the higher the density gradient. Although the model itself received a lot of critique, it is nevertheless interesting to look at the structure of the associated system.  The Brenner model, like the compressible Euler and the Navier-Stokes equations, consists of two conservation equations
\begin{subnumcases}{\label{Brenner}}
\partial_t \vr + \Div (\vr\, \vv_m) = 0,\label{B1}\\
\partial_t (\vr \,\vv) + \Div (\vr \,\vv \otimes \vv_m ) -\Div \mathbb{T}(\vr,\vv)= 0,\label{B2}
\end{subnumcases} 
with the closure relation $\vv_m=\vv-K\Grad \log\vr$. Above $\mathbb{T}$ is some general stress tensor that might include the presssure component, and  we do not wish to specify it here. In the dissipative AR system the closure relation between the velocities is simply $\vv_m=\vv-\Grad p(\vr)$. Surprisingly, quite a few other systems considered in the literature could be put in the same framework, including:
\begin{itemize}
    \item compressible Navier-Stokes with density-dependent viscosity $\mu(\vr)$, for which $\vv_m=\vv-\frac{ \mu'(\vr)}{\vr}\Grad\vr$ investigated, for example, in \cite{BDL, VY, VY2, BVY};
    \item compressible Navier-Stokes equations with transport of noise, for which $\vv_m=\vv-\mathbb{Q}\circ {\rm{d}} \mathbf{W}$, see \cite{Holm,FL,BFHZ};
    \item odd viscosity model with $\vc{v}_m=\vv-2\nu_0\Grad^\perp \log\vr$, see \cite{FaVa}.
\end{itemize}
Moreover, if the local offset function $\Grad p(\vr)$ in dissipative AR system \eqref{AR_multi},  is replaced by the non-local one $\Grad K \ast  \vr$, then the system becomes the Euler-alignment model with matrix-valued communication weight. As explained in \cite{CPSZ}, for sufficiently regular solutions, and sufficiently smooth and symmetric kernels $K(x)$, the dissipative AR system is equivalent to
\begin{subnumcases}{\label{EA}}
   	\pt \vr+\Div (\vr \vu)=0, \label{EA1} \\
\pt(\vr \vu )+ \Div(\vr \vu \otimes \vu)- \vr  \lr{\int_{\R^d} \lr{\vu(y)- \vu(x)}^{\intercal}\boldsymbol{\Psi}(x-y)  \vr(y)\, \dy } =0, \label{EA2}
\end{subnumcases}
where 
\eqh{\label{khess}
\boldsymbol{\Psi}=\nabla^2 K = \lr{\partial^2_{x_i x_j}K }_{i,j=1}^d.}

The equivalence between the nonlocal dissipative AR models and the Euler-alignment system \eqref{EA} has been investigated also at the levels of the corresponding particle and kinetic formulations. In simplified variants -- in one space dimension or for non-singular kernels -- it has been used  to predict behaviour of the solutions, for example in \cite{HPZ, HKPZ,Kim21,CZ, Kim}. In the case of singular communication weight, corresponding roughly to  $K(x)=|x|^{2-\alpha}$,  the equivalence was obtained by Peszek and Poyato in \cite{PP22}.

In the last part of this section we discuss the connection between the one-dimensional dissipative AR model and the lubrication model. The effective continuous model of lubrication was derived from the particle description by Lefebvre-Lepot and  Maury in \cite{LM}. In case when inertia of the flow is taken into account the derived equations are of the form:
\begin{subnumcases}{\label{NS_1Dsys}}
\partial_t \vr + \partial_x(\vr u) = 0, \label{1D1} \\
		\partial_t(\vr u)+\partial_x(\vr u^2)-\partial_x\lr{\lambda_\ep(\vr)\partial_x u}=0, \label{1D2}
\end{subnumcases}
where \[\lambda_\ep(\vr)\approx \frac{\ep}{1-\vr}\] is the effective viscosity. It is proportional to the viscosity of the interstitial fluid $\ep$, and it takes into account the lubrication forces, i.e. forces preventing contact between grains: $\lambda(\vr)\to +\infty$ as $\vr\to 1$, where $1$ corresponds to fully congested flow. Comparing  equations \eqref{1D2}and  \eqref{NS_1D}, it is clear that the lubrication model is a special one-dimensional case of the dissipative AR model for the particular choice  of the cost function  $$p(\vr)=\int_{0}^\vr\frac{\lambda_\ep(s)}{s^2} {\rm d}s.$$ 
In the articles \cite{CNPZ,CMPZ} our purpose was to investigate the global in time existence of solutions to system \eqref{NS_1Dsys} on one-dimensional torus, and on the whole line, and to justify the hard-congestion limit that corresponds to $\ep\to0$. Another way to obtain the hard-congestion limit in the model is to consider $\lambda_n(\vr)=\vr^{\gamma_n}$ with $\gamma_n\to\infty$ for $n\to\infty$, which approximates the barrier $\vr=1$ from above, and not from below like the $\ep$-approximation, see \cite{M} for a relevant result.

For particular approximation of the viscosity coefficient $\lambda_\ep$, i.e. 
$$\lambda_\ep(\vr_\ep)= \vr_\ep^2 p'_\ep(\vr_\ep)+\vr_\ep^2 \varphi'_\ep(\vr_\ep),$$
where $$
p_\ep(\vr_\ep) = \ep\dfrac{\vr_\ep^\gamma}{(1-\vr_\ep)^\beta}, \quad 
\varphi_\ep(\vr_\ep)  = \dfrac{\ep}{\alpha-1} \vr_\ep^{\alpha-1}, \quad
\gamma \geq 0, \quad \beta > 1, \quad \alpha \in \lr{0,\frac{1}{2}},
$$
we showed  in \cite{CNPZ} that when $\ep>0$ is fixed, the strong solutions exist, are unique and global-in-time. Moreover, it was shown that $0 < \vr_\ep(t,x) \leq C_\ep<1$ with $C_\ep\to1$ as $\ep\to 0$. The proof followed the earlier results of Burtea and Haspot \cite{BH} and of Constantin {\em et al.} \cite{Constantin_2020}. The problem with passing to the limit $\ep\to0$ lies  precisely in lack of the uniform (with respect to $\ep$) bound for the density away from 1, which causes problems in passage to the limit in the diffusion term. This, together with lack of estimates for $u_\ep$  independent of $\vr_\ep$ makes the limit in the convective term $\vr_\ep u_\ep^2$ difficult to identify. The same happens for the reformulation of the system in the dissipative AR form: there is no estimate  for $w_\ep=u_\ep+\px p_\ep(\vr_\ep)$, and the limit of the convective term $\vr_\ep u_\ep w_\ep$ might be some general measure. Formally, the hard-congestion limit was identified already in \cite{LM}, and it agrees with one-dimensional reduction of the limit from the recent result of Aceves et al. \cite{ABDM}. Denoting by $\pi'_\ep(\vr_\ep)=\vr_\ep p'(\vr_\ep)$, one expects the following asymptotic behaviour:
\eq{\label{formal}\left\{
\begin{array}{l}
\partial_t \vr_\ep + \partial_x(\vr_\ep u_\ep) = 0  \\
		\partial_t (\vr_\ep w_\ep) + \partial_x (\vr_\ep w_\ep u_\ep) = 0\\ 
		{\vr_\ep w_\ep} = {\vr_\ep u_\ep} + \px \pi_\ep(\vr_\ep)  
\end{array}\right.
\quad 
\begin{array}{c}
{\longrightarrow}\\
{\ep\to0}
\end{array}
\quad
\left\{
\begin{array}{l}
\partial_t \vr + \partial_x(\vr u ) = 0  \\
		\partial_t (\vr w ) + \partial_x (\vr w  u ) = 0\\ 
		{\vr  w } = {\vr  u } + \partial_x \pi \\
		0\leq \vr\leq1,\ \pi\geq 0,\ (1-\vr)\pi=0,
\end{array}\right.
}
where the limit of $\pi_\ep$, denoted by $\pi$, is a new unknown in the system on the right, satisfying the algebraic constraint $(1-\vr)\pi=0$. This constraint means that the support of  $\pi$ is contained in the set where $\vr=1$, i.e. where the congestion appears. Physically, it corresponds to appearance of an additional forcing term preventing overcrowding.

To justify this limit passage rigorously, we investigated two different regimes of global-in-time solutions:  weak and measure solutions. In either of the cases, it is essential to obtain some further estimates for the velocity $u_\ep$, independently  of the a-priori bounds for the $\ep$-modulated viscosity. Do this end, we first note that the maximum principle for the \emph{active potential} $V_\ep=\lambda_\ep(\vr_\ep)\partial_x u_\ep$, satisfying formally 
$$\pt V_\ep+a_\ep \px V_\ep-b_\ep\px^2 V_\ep=-c_\ep V^2_\ep,$$
with some coefficients  $b_\ep,\ c_\ep\geq0$, implies the
one-sided Lipschitz condition  on $u_\ep$:
$${\px u_\ep\leq C.}$$
Although, this condition is still insufficient to pass to the limit in the term $\vr_\ep u_\ep (u_\ep+\px \pi_\ep)$, one can proceed further in two different ways:
\begin{itemize}
\item In the finite-energy approach \cite{CNPZ}, by essentially integrating the momentum equation, one can derive additional uniform in $\ep$ bound for $\pi_\ep$ in $L^\infty(0,T; H^1(\T))$. This, combined with standard compensated compactness argument, allows to pass to the limit in the troublesome term in the sense of distributions;

\item In  the duality formulation \cite{CMPZ}, introduced initially by Bouchut and James \cite{BJ98, BJ99} for pressureless gases dynamics equations, and extended by Boudin to study  the viscous approximation thereof in \cite{Boudin}, one solves the same problem by the appropriate choice of the test function which makes the convective term disappear from the weak formulation. The adequate choice is the {\emph{reversible}} solution of the dual transport equation with the velocity vector field $u_\ep$, which satisfies the Oleinik entropy condition and the uniform $L^\infty$ bound. These, in turn,  guarantee the uniqueness of solution.  Justification of limit \eqref{formal}  is therefore reduced to showing the stability of the duality solutions.
\end{itemize}


\bibliographystyle{plain}

\bibliography{biblio}

\end{document}